\documentclass{article}
\usepackage{amsmath}
\usepackage{amssymb,comment,amscd}
\usepackage{mathrsfs,wasysym}

\usepackage{euscript,color}

\usepackage{libertine}
\usepackage[T1]{fontenc}

\usepackage[utf8]{inputenc}

\usepackage{amsmath}
\usepackage{amssymb}
\usepackage{amsthm}
\usepackage{mathrsfs}

\usepackage{graphicx}
\graphicspath{ {/desktop/graphs/} }

\usepackage{version}
\usepackage[margin=3cm]{geometry}
\usepackage{color}
\usepackage{mathtools}
\usepackage{tikz}
\usetikzlibrary{arrows}
\usetikzlibrary{trees}
\usepackage{float}
\usepackage{textcomp}
\usepackage{bbm}

\numberwithin{equation}{section}

\newtheorem{proposition}{Proposition}[section]
\newtheorem{claim*}{Claim}[section]

\newtheorem{definition}{Definition}[section]
\newtheorem{lemma}{Lemma}[section]
\newtheorem{theorem}{Theorem}[section]

\newtheorem{assumption}{Assumption}

\newtheorem{remark}{Remark}

\newcommand{\rr}{\mathbb R}

\newcommand{\bs}{\boldsymbol}
\usepackage{bm}
\newcommand{\vt}{\vartheta}

\newcommand{\vb}{\vspace{2mm}}
\newcommand{\s}{^\star}

\newcommand*\longbar[1]{%
  \hbox{%
    \vbox{%
      \hrule height 0.4pt % The actual bar
      \kern0.4ex%         % Distance between bar and symbol
      \hbox{%
       \kern-0.04em%      % Shortening on the left side
        \ensuremath{#1}%
        \kern-0.04em%      % Shortening on the right side
      }%
    }%
  }%
}

\renewcommand{\check}{\widetilde}
\renewcommand{\overline}{\longbar}
\renewcommand{\bar}{\overline}

\begin{document}

\title{Dynamic Erd\H{o}s-R\'enyi graphs}
\author{M.~Mandjes, N.~J.~Starreveld,  R.~ Bekker and P.~Spreij}
\maketitle

\begin{abstract}
We propose two classes of dynamic versions of the classical Erd\H{o}s-R\'enyi graph: one in which the transition rates are governed by an external regime process, and one in which the transition rates are periodically resampled. For both models we consider the evolution of the number of edges present, with explicit results for the corresponding moments, functional central limit theorems and large deviations asymptotics. 
\end{abstract}

\section{Introduction}
Over the past decades, networks have been the subject of an intensive research effort. As networks offer the right framework to model e.g.\ social, physical, chemical, biological and technological phenomena, various specific aspects have been studied in depth. Arguably among the most studied objects is the {\it Erd\H{o}s-R\'enyi graph} \cite{ER,EG}. In such a random graph $G(n,p)$ there are $n$ vertices, and each of the $N=\binom{n}{2}$ edges is `up' with a fixed probability $p\in(0,1)$ or `down' otherwise. By now there is a sizeable literature on this type of graph, providing detailed insight into its probabilistic properties, an example of a key result being that if the `up-probability' $p$ is larger than $\log n/n$, then the resulting graph is almost surely connected.

The existing  literature predominantly focuses on {\it static} graphs: the random graph is drawn just once, and does not change over time. In many real-life situations, however, the network structure temporally evolves, with edges appearing and disappearing. In a few recent contributions, first results on such dynamic random graphs have been reported, but the analysis of this class of models is still in its infancy; see e.g.\ \cite{HO1,HO2,ZMN}, and \cite{BBRJ} for an illustration of its use in engineering.  

In \cite{ZMN} various dynamic random graph models are discussed, among them a dynamic Erd\H{o}s-R\'enyi graph in which all $N$ edges evolve independently. In this model, each edge makes transitions from present to absent and vice versa in a Markovian manner: it exists for an exponential time with parameter $\mu$ (which we refer to as the `up-rate'), and disappears for an exponential time with parameter $\lambda$ (the `down-rate'). For this model various metrics can be analyzed in closed form. In particular the distribution of the number of edges at time $t$, throughout this paper denoted by $Y(t)$, can be explicitly computed. A special case is that in  which no edges exist at $t=0$: then the distribution of $Y(t)$  coincides with the number of edges in a static Erd\H{o}s-R\'enyi graph $G(n,p(t))$ (with an up-probability that depends on $t$).

\vb

In many applications the model that we just sketched is of limited relevance, as various features that play a role in real-life networks are not covered.  To remedy this, in \cite{ZMN} alternative random graph processes were proposed, such as the dynamic counterparts of the configuration model and the stochastic block model.  It is noted that a specific property that is often not fulfilled in real networks is that of the edges evolving independently; in practice likely there will be `external' factors that affect all these $N$ processes simultaneously, rendering them dependent. An example is a dynamic random graph in which the values of the up-rate and down-rate are determined by an independent stochastic process (think of temperature in a chemical network, weather conditions in a road traffic network, economic conditions in a financial network, etc.). 

Motivated by the above considerations, the focus of this paper is on models in which the edges evolve {\it dependently}; the main contribution is 
that we propose and analyze two such models. In the first model, studied in Section \ref{RS}, the up-rate and the down-rate of each of the edges are determined by an external, autonomously evolving Markov process $X(t)$, in the sense that at time $t$ these rates (for all edges) are $\lambda_i$ and $\mu_i$ if $X(t)=i$; this mechanism is usually referred to as {\it regime switching}. In the second model, which is analyzed in Section \ref{ERR}, the up-rate and the down-rate (say, $\Lambda$ and $M$) are resampled every $\Delta>0$ time units (and these sampled values then apply to all edges). 

In more detail, our findings are the following. The focus is on the probabilistic properties of the process $Y(t)$ that records the number of edges present as a function of time. For both models mentioned above we manage to uniquely characterize its transient and stationary behavior, albeit in a somewhat implicit way: for the first model in terms of a {\sc pde} for the corresponding probability generating function ({\sc pgf}), for the second model in terms of a recursion for the {\sc pgf}. Then we use these characterizations to point out how transient and stationary means can be computed. The next step is to consider scaling limits; under a particular scaling, the process $Y(t)$ satisfies a functional central limit theorem. More specifically, after centering and scaling it converges to an Ornstein-Uhlenbeck ({\sc ou}) process; interestingly, in \cite{TYA} it is shown that for certain dynamic Erd\H{o}s-R\'enyi graphs that a particular clique-complex related quantity (the `Betti number') is described by an {\sc ou} process as well. Finally we discuss for both models the corresponding sample-path large deviations, characterizing the models' rare-event behavior.
In Section \ref{NUM}, the results are illustrated by numerical examples.

\section{Erd\H{o}s-R\'enyi graphs under regime switching} \label{RS}
In this section we consider the following model. Let $(X(t))_{t\geqslant 0}$ be an irreducible continuous-time Markov process, typically referred to as the {\it regime process} or {\it background process}, living on the state space $\{1,\ldots,d\}$. The transition rate matrix corresponding to $(X(t))_{t\geqslant 0}$ is denoted by $Q=(q_{ij})_{i,j=1}^d$ and the corresponding invariant distribution by the (column) vector ${\boldsymbol \pi}$. As before, we consider the situation of $N$ possible vertices. Let $\mu_i\geqslant 0$ be the hazard rate of an existing edge becoming inactive when the regime process is in state $i$; likewise, $\lambda_i\geqslant 0$ is the hazard rate corresponding with a non-existing edge becoming active. Due to the common regime process the edges are reacting to, the number of links present (denoted by $(Y(t))_{t\geqslant 0}$) evolves according to an interesting dynamic structure. 

\subsection{Generating function}
We start our exposition by studying the (transient and stationary) {\sc pgf}\,s  
\[\phi_i(t,z):={\mathbb E} \left(z^{Y(t)}\,1_{\{X(t) = i\}}\right),\:\:
\phi_i(z):={\mathbb E} \left(z^{Y}\,1_{\{X= i\}}\right).\]
We do so by first analyzing $p_i(m, t):={\mathbb P}(Y(t) = m, X(t) = i),$
by following classical procedures; later we also point out how $p_i(m):={\mathbb P}(Y = m, X = i)$ can be found.
Setting up the Kolmogorov equations, with $q_i:=-q_{ii}>0$,
\begin{eqnarray*} 
p_i(m,t+\Delta t)&=& \sum_{j\not=i} p_j(m,t) q_{ji}\,\Delta t \\&&+\, p_i(m+1,t) \mu_i (m+1)\Delta t
+ p_i(m-1,t) \lambda_i(N-m+1)\Delta t \\&& +\, p_i(m,t)\big(1- q_i\Delta t - \mu_i\,m\,\Delta t-\lambda_i\,(N-m)\,\Delta t\big)+o(\Delta t),
\end{eqnarray*}
leading to the linear system of differential equations 
\begin{eqnarray*}p_i'(m,t) &=& \sum_{j=1}^d p_j(m,t)q_{ji} + p_i(m+1,t) \mu_i\,(m+1)
 \\
&&+\,p_i(m-1,t) \lambda_i\,(N-m+1) -p_i(m,t)\mu_i\,m-p_i(m,t)\lambda_i\,(N-m),
\end{eqnarray*}
where $p_i(-1,t)$ and $p_i(N+1,t)$ are set to $0$.
Multiplying by $z^m$ and summing over $m=0$ up to $N$, we arrive at the {\sc pde}
\begin{eqnarray*}
\frac{\partial}{{\partial t}}\phi_i(t,z) &=& \sum_{j=1}^d \phi_j(t,z)q_{ji}+\mu_i (1-z)\frac{\partial}{{\partial z}}\phi_i(t,z)+\\&&\lambda_iN(z-1) \phi_i(t,z)+\lambda_iz(1-z)\frac{\partial}{{\partial z}}\phi_i(t,z).
\end{eqnarray*}
In stationarity, the left-hand side of the previous display can be equated to 0, thus leading to an {\sc ode}. We obtain
\[0=\sum_{j=1}^d \phi_j(z)q_{ji}+\mu_i (1-z)\phi_i'(z)+\lambda_iN(z-1) \phi_i(z)+\lambda_iz(1-z)\phi_i'(z).\]

\subsection{Moments}
Following a standard procedure, we can find explicit expressions for all (factorial) moments. To this end, we define $e_{i,k}:= {\mathbb E}((Y)_k1_{\{X=i\}})$, with $(x)_k$ denoting $x(x-1)\cdots(x-k+1)$. We obtain the factorial moments  by
differentiating with respect to $z$ and plugging in $z=1$: in self-evident matrix/vector notation, with
$\Lambda:={\rm diag}\{{\boldsymbol \lambda}\}$ and $M:={\rm diag}\{{\boldsymbol \mu}\}$,
\[{\boldsymbol 0}^{\rm T} = {\boldsymbol e}_1^{\rm T}Q - {\boldsymbol e}_1^{\rm T}M +{\boldsymbol \pi}^{\rm T}\Lambda N
- {\boldsymbol e}_1^{\rm T}\Lambda.\]
This leads to ${\mathbb E}Y =  {\boldsymbol e}_1^{\rm T}{\boldsymbol 1}$, with 
${\boldsymbol e}_1^{\rm T}= N\cdot{\boldsymbol \pi}^{\rm T}\Lambda (\Lambda+M-Q)^{-1};$
observe that the mean is proportional to $N$, as expected. This procedure can be found to find a recursion for 
all factorial moments: by differentiating $k$ times and inserting $z=1$, we obtain, for $k=2,3,\ldots,N$,
\[{\boldsymbol 0}^{\rm T} = {\boldsymbol e}_k^{\rm T}Q - k\,{\boldsymbol e}_k^{\rm T}M +kN\,{\boldsymbol e}_{k-1}^{\rm T}\Lambda 
-k\, {\boldsymbol e}_k^{\rm T}\Lambda -  k(k-1)\,{\boldsymbol e}_{k-1}^{\rm T}\Lambda,\]
and consequently
\[{\boldsymbol e}_k^{\rm T} = k\,(N-k+1)\cdot{\boldsymbol e}_{k-1}^{\rm T} \,\Lambda (k\Lambda+kM-Q)^{-1}.\]
Observe that this recursion can be explicitly solved, as we know ${\boldsymbol e}_1^{\rm T}$; the following result now straightforwardly follows.
\begin{proposition} For $k=1,\ldots,N$, 
\[{\boldsymbol e}_k^{\rm T} = k!\,(N)_k \cdot {\boldsymbol \pi}^{\rm T}  \Lambda (\Lambda+M-Q)^{-1}\Lambda (2\Lambda+2M-Q)^{-1}\cdots \Lambda (k\Lambda+kM-Q)^{-1},\]
whereas ${\boldsymbol e}_k^{\rm T} =0$ for $k=N+1,N+2,\ldots$. 
\end{proposition}

Following standard techniques, we can now evaluate all stationary probabilities as well. First, $p_i(N)$ follows from the identity
${e}_{i,N} = {\mathbb E}((Y)_N1_{\{X=i\}})= N! \,p_i(N).$
We can recursively find the other probabilities $p_i(m)$; applying
\[e_{i,N-1} =  {\mathbb E}((Y)_{N-1}1_{\{X=i\}})= (N-1)! \, p_i(N-1) + N!\, p_i(N),\]
we can express $p_i(N-1)$ in terms of $p_i(N)$ (and ${e}_{i,N-1}$  and ${e}_{i,N}$). In general $p_i(m)$ can be found from
$p_i(m+1),\ldots,p_i(N)$ using
\[e_{i,m} = \sum_{k=m}^N (k)_m p_i(k).\]

\begin{remark}
In addition, the {\it transient} factorial moments ${\mathbb E}((Y(t))_k\,1_{\{X(t)=i\}})$ can be (recursively) found;
in every step of the recursion a system of linear differential equations (rather than a linear-algebraic equation) needs to be solved; see \cite{MT} for a similar procedure
in the context of infinite-server queues under regime switching.\end{remark}

\subsection{Diffusion results under scaling}
In this subsection we impose the scaling $Q\mapsto N^\delta Q$, entailing that the regime process is sped up by a factor $N^\delta$, with the objective to prove a functional central limit theorem for the resulting limiting process. To get a feeling for how this scaling affects the system's behavior, we first compute the mean and variance of the stationary number of edges. To this end, we use the following lemma, which is proven in the appendix. In the sequel $D:=({\bs 1}{\bs \pi}^T-Q)^{-1}-{\bs 1}{\bs \pi}^T$ denotes the {\it deviation matrix}. Also
$x\s:= {\boldsymbol x}^{\rm T}{\boldsymbol\pi}$ for ${\bs x}\in{\mathbb R}^d$ and $\Gamma:={\rm diag}\{{\boldsymbol \gamma}\}=\Lambda+M$. Let ${\bs\gamma}:={\bs\lambda}+{\bs\mu}$ be componentwise positive.

\begin{lemma} \label{lalemma} 
Define $F_{N,k}:= (k\,\Gamma -NQ)^{-1}$ for $k\in{\mathbb N}$. Then, as $N\to\infty$,
\[F_{N,k} = \frac{1}{k}\frac{1}{\gamma\s} {\boldsymbol 1}\,{\boldsymbol\pi}^{\rm T}+
\frac{1}{N}E+O(N^{-2}),\:\:E:=\left(I-\frac{1}{\gamma\s}{\bs 1}\,{\bs \pi}^{\rm T}\Gamma\right) D\left(I- \frac{1}{\gamma\s}
{\bs\gamma}\,{\bs\pi}^{\rm T}\right)\hspace{-0.5mm}.\]
\end{lemma}
Let us first evaluate the mean of $Y$ under this scaling; in the steps below we use ${\bs \pi}^{\rm T}\Lambda\, {\bs 1}=\lambda\s$ and $D{\bs 1}={\bs 0}$. From the above lemma, we find, with  $\bar\varrho:=\lambda\s/\gamma\s$,
\begin{eqnarray*}
{\mathbb E}\,Y&=&N {\boldsymbol \pi}^{\rm T}  \Lambda (\Lambda+M-N^\delta Q)^{-1}{\bs 1}
= N {\boldsymbol \pi}^{\rm T}  \Lambda  F_{N^\delta,1}{\bs 1}
\\&=& N {\boldsymbol \pi}^{\rm T}  \Lambda\left(\frac{1}{\gamma\s} {\boldsymbol 1}{\boldsymbol\pi}^{\rm T}{\bs 1}+
{N}^{-\delta}\,E{\bs 1}+O(N^{-2\delta})\right)=N\,\bar\varrho+
%N^{1-\delta} {\boldsymbol \pi}^{\rm T}  \Lambda \,E{\bs 1}+
%o(N^{1-\delta}) .\end{eqnarray*}
O(N^{1-\delta}) .\end{eqnarray*}
Along the same lines,
\[({\mathbb E}Y)^2 = N^2\bar\varrho\,^2 - N^{2-\delta}\frac{2}{\gamma\s} {\bs \pi}^{\rm T}  (\Lambda-\bar\varrho\,\Gamma) D\,  \bar\varrho \,\Gamma{\bs 1} +o(N^{{\rm max}\{1,2-\delta\}}).\]
In addition, ignoring sublinear terms,
\begin{eqnarray*}
{\mathbb E}Y(Y-1) &=&2N(N-1)\,{\boldsymbol \pi}^{\rm T}  \Lambda (\Lambda+M-N^\delta Q)^{-1}\Lambda (2\Lambda+2M-N^\delta Q)^{-1}{\bs 1}\\
&=& 2N(N-1)\,{\boldsymbol \pi}^{\rm T}  \Lambda\, F_{N^\delta,1} \, \Lambda \,F_{N^\delta,2}{\bs 1}\\
&=&2N(N-1) \,{\boldsymbol \pi}^{\rm T}  \Lambda\, \left(\frac{1}{\gamma\s} {\boldsymbol 1}\,{\boldsymbol\pi}^{\rm T}+
\frac{1}{N^\delta}E\right)\Lambda\left(\frac{1}{2\gamma\s} {\boldsymbol 1}\,{\boldsymbol\pi}^{\rm T}+
\frac{1}{N^\delta}E\right){\bs 1}.
\end{eqnarray*} 
Using the following equalities 
\begin{eqnarray*}
{\boldsymbol \pi}^{\rm T}  \Lambda\, \left(\frac{1}{\gamma\s} {\boldsymbol 1}\,{\boldsymbol\pi}^{\rm T}\right)\Lambda\left(\frac{1}{2\gamma\s} {\boldsymbol 1}\,{\boldsymbol\pi}^{\rm T}\right){\bs 1} &=&\frac{\bar\varrho\,^2}{2},\\
{\boldsymbol \pi}^{\rm T}  \Lambda\, E\Lambda \left(\frac{1}{2\gamma\s} {\boldsymbol 1}\,{\boldsymbol\pi}^{\rm T}\right){\bs 1} &=&\frac{1}{2\gamma\s} {\bs \pi}^{\rm T} (\Lambda-\bar\varrho\,\Gamma) D (\Lambda-\bar\varrho\,\Gamma) {\bs 1},\\
{\boldsymbol \pi}^{\rm T}  \Lambda\, \left(\frac{1}{\gamma\s} {\boldsymbol 1}\,{\boldsymbol\pi}^{\rm T}\right)\Lambda E{\bs 1}&=&-\frac{1}{\gamma\s} {\bs \pi}^{\rm T}  (\Lambda-\bar\varrho\,\Gamma) D \, \bar\varrho \,\Gamma{\bs 1},
\end{eqnarray*}
we arrive at
\[{\mathbb E}Y(Y-1) =N(N-1)\,\bar\varrho\,^2+N^{2-\delta}\frac{1}{\gamma\s} {\bs \pi}^{\rm T} (\Lambda-\bar\varrho\,\Gamma) D (\Lambda-3\,\bar\varrho\,\Gamma) {\bs 1}+o(N^{{\rm max}\{1,2-\delta\}}).\]
By virtue of the identity ${\mathbb V}{\rm ar} \,Y= {\mathbb E}Y(Y-1) +{\mathbb E}Y-({\mathbb E}Y)^2$, we thus find
\begin{equation}
\label{VARi}{\mathbb V}{\rm ar} \,Y = N\,\bar\varrho(1-\bar\varrho) + N^{2-\delta}\,v\,+o(N^{{\rm max}\{1,2-\delta\}}),\end{equation}
{with}\[v:=\frac{1}{\gamma\s} {\bs \pi}^{\rm T} (\Lambda-\bar\varrho\,\Gamma) D (\Lambda-\bar\varrho\,\Gamma) {\bs 1}.\]
It can be checked that this formula is symmetric, in the sense that it is invariant under swapping ${\bs \lambda}$ and ${\bs \mu}$, which is in line with ${\mathbb V}{\rm ar}\,Y={\mathbb V}{\rm ar}\,(N-Y)$; note that $\Lambda-\bar\varrho\,\Gamma =
(1-\bar\varrho)\Lambda-\bar\varrho M$.

Upon inspecting the asymptotic shape of ${\mathbb V}{\rm ar}\,Y$, we observe a dichotomy. For $\delta>1$ the regime process jumps so fast that all edges essentially behave independently, experiencing an `effective up-rate' of $\lambda\s$, and  an `effective down-rate' of $\mu\s$, so that in this regime $Y$ is approximated with a Binomial random variable with parameters $N$ and $\bar\varrho$. For $\delta<1$ the regime process is relatively slow, and hence affects the variance (which is, as a result, superlinear in $N$). 

\vb

We now prove a functional central limit theorem. For the moment we focus on the case $\delta=1$; in Remark \ref{rem2} we comment on what happens when $\delta>1$ or $\delta<1$. Let $P_1(\cdot)$ and $P_2(\cdot)$ be two independent unit-rate Poisson processes. With $Z_i(s):=1_{\{X(s)=i\}}$, and $Y(0)=0$ (remarking that any other starting point can be dealt with similarly), 
\begin{equation}\label{YT}Y(t) =   P_1\left(\sum_{i=1}^d \int_0^t \lambda_i Z_i(s)(N- Y(s)){\rm d}s\right) -P_2\left(
\sum_{i=1}^d \int_0^t \mu_iZ_i(s) Y(s){\rm d}s\right) .\end{equation}
The first step is to verify that $Y(t)/N$ converges to $y(t)$, defined as the solution of the integral equation
\[y(t) =  \lambda\s  \int_0^t (1-y(s)){\rm d}s -\mu\s \int_0^t y(s){\rm d}s,\]
i.e., $y(t) =   \varrho(t):=\bar\varrho\cdot (1-e^{-\gamma\s t}).$ Define
\begin{equation}\label{barY}\bar Y(t):= \frac{Y(t) - N\varrho(t)}{\sqrt{N}};\end{equation}
our objective is to prove that $\bar Y(\cdot)$ converges to a Gaussian process (and we identify this process). As we follow \cite[Section 5]{BTM},
which in turn uses intermediate results of \cite{HJMST}, we restrict ourselves to the most important steps.

We know from (\ref{YT}) that, for some martingale $K(t)$,
\[{\rm d}Y(t) = {\bs\lambda}^T {\bs Z}(t) (N-Y(t)) {\rm d}t - {\bs\mu}^{\rm T} {\bs Z}(t) Y(t) {\rm d}t +{\rm d}K(t),\]
and therefore
\[{\rm d}\bar Y(t) = \sqrt{N}\big((1-\varrho(t)){\bs\lambda}^{\rm T}-\varrho(t){\bs\mu}^{\rm T}\big){\bs Z}(t)  {\rm d}t
-{\bs\gamma}^{\rm T}{\bs Z}(t) \bar Y(t){\rm d}t +\frac{{\rm d}K(t)}{\sqrt{N}} -\sqrt{N}\varrho'(t){\rm d}t.\]
Now define
$W(t):=e^{Z_+(t)}\bar Y(t),$ where ${Z}_+(t):=\int_0^t 
{\bs \gamma}^{\rm T}{\bs Z}(s){\rm d}s,$
so that, 
\[{\rm d}W(t) =e^{Z_+(t)}\left(\sqrt{N}\big((1-\varrho(t)){\bs\lambda}^{\rm T}-\varrho(t){\bs\mu}^{\rm T}\big){\bs Z}(t)  {\rm d}t
+\frac{{\rm d}K(t)}{\sqrt{N}} -\sqrt{N}\varrho'(t){\rm d}t
\right) .\]
Observing that 
$\big((1-\varrho(t)){\bs\lambda}^{\rm T}-\varrho(t){\bs\mu}^{\rm T}\big){\bs \pi}= \varrho'(t),$
and recalling that ${\bs\gamma}={\bs\lambda}+{\bs \mu}$, the equality in the previous display simplifies to
\[{\rm d}W(t) = e^{Z_+(t)}\left(\sqrt{N}\big({\bs\lambda}^{\rm T}-\varrho(t){\bs\gamma}^{\rm T}\big)({\bs Z}(t) -{\bs\pi}) {\rm d}t
+\frac{{\rm d}K(t) }{\sqrt{N}}
\right) .\]
We now consider the two terms in the previous display separately. As was established in \cite{BTM,HJMST}, for the first term, as $N\to\infty$,
\[\int_0^\cdot \sqrt{N} e^{Z_+(s)}
\big({\bs\lambda}^{\rm T}-\varrho(s){\bs\gamma}^{\rm T}\big)({\bs Z}(s) -{\bs\pi}) {\rm d}s\to \int_0^\cdot e^{\gamma \s s}{\rm d} G(s),\]
where $G(\cdot)$ satisfies
\begin{equation}\label{defg1}\langle G\rangle_t = g(t):=2\int_0^t {\bs\pi}^{\rm T} (\Lambda -\varrho(s)\Gamma) D  (\Lambda -\varrho(s)\Gamma) {\bs 1}{\rm d}s.\end{equation}
Also as in \cite{BTM,HJMST}, the second term obeys, as $N\to\infty$,
\[\int_0^\cdot \frac{1}{\sqrt{N}} e^{Z_+(s)} {\rm d}K(s)  \to 
 \int_0^\cdot e^{\gamma\s s}{\rm d} H(s),\]
where $H(\cdot)$ satisfies (using the relation between $K(\cdot)$ and the Poisson processes $P_1(\cdot)$ and $P_2(\cdot)$)
\begin{equation}\label{defh1}\langle H\rangle_t = h(t):= \int_0^t \lambda\s(1-\varrho(s)){\rm d}s + \int_0^t \mu\s \varrho(s){\rm d}s.\end{equation}
Combining the two terms studied above, it thus follows that, as $N\to\infty$, $W(\cdot)$ weakly converges to $W_\infty(\cdot)$, which is the solution to the stochastic differential equation, with $B(\cdot)$ a standard Brownian motion, 
\begin{equation}
\label{defW}{\rm d}W_\infty(t) =e^{\gamma\s t} \sqrt{g'(t)+h'(t)} \,{\rm d}B(t).\end{equation}
Translating this back  in terms of a stochastic differential equation, again mimicking the line of reasoning of \cite{BTM,HJMST}, we obtain the following result.

\begin{theorem}
$\bar Y(\cdot)$ converges weakly to $\bar Y_\infty(\cdot)$, which is the solution to the stochastic differential equation
\begin{equation}\label{limde}{\rm d}\bar Y_\infty(t)= - \gamma\s\, \bar Y_\infty(t)\,{\rm d}t + \sqrt{g'(t)+h'(t)} \,{\rm d}B(t),\end{equation}
with  
$g(\cdot)$ and $h(\cdot)$ given by $(\ref{defg1})$ and $(\ref{defh1})$, respectively.
\end{theorem}

\begin{remark}
Using the behavior of $g'(t)$ and $h'(t)$ for $t$ large, we conclude that for large values of $t$ (`in stationarity'), this stochastic differential equation reads
\[ {\rm d}\bar Y_\infty(t)= - \gamma\s\, \bar Y_\infty(t)\,{\rm d}t + \sqrt{2\gamma\s\,\bar\varrho\,(1-\bar\varrho)+2\gamma\s\,v}\,{\rm d}B(t),\]
which defines an {\sc ou} process with mean $0$ and variance $\bar\varrho\,(1-\bar\varrho) + v$; note that this aligns with what we found, plugging in $\delta=1$, in (\ref{VARi}).\end{remark}

\begin{remark} \label{rem2} When $\delta<1$, the $\sqrt{N}$ in the definition of (\ref{YT}) needs to be replaced by $N^{\delta/2}$; it is readily checked that in the limiting stochastic differential equation~(\ref{limde}) we then just have $g'(t)$ below the square-root sign. On the contrary, if $\delta>1$ then the definition of~(\ref{YT}) remains unchanged, but below the square-root sign in (\ref{limde}) we only have $h'(t)$.
\end{remark}

\subsection{Large deviations results under scaling} Where we above discussed the diffusion behavior of the process under study, we now consider rare events. We again focus on the scaling corresponding to $\delta=1$, following the setup of \cite{HMS}. Intuitively, the rare-event behavior is decomposed into the effect of the regime process, and that of the edge dynamics conditional on the regime process. 

Let ${\bs g}(\cdot)$ be in  $U_T$, defined as  the set of non-negative $d$-dimensional functions such that the $g_i(s)$  sum to $1$, for all $s\in[0,T]$. Then
\[{\mathbb J}_T({\bs g}) := \int_0^T \sup_{{\bs u}\geqslant {\bs 0}}\left(-\sum_{i=1}^d \frac{(Q{\bs u})_i}{u_i}g_i(s)\right){\rm d}s.\]
In addition,
\[\Lambda_{x,{\bs g}}(\vt):=\sum_{i=1}^d g_i\left(x\lambda_i(e^\vt-1)+(1-x)\mu_i(e^{-\vt}-1)\right).\]
Based on the findings in \cite{HMS}, one anticipates a sample-path {\sc ldp} (of `Mogulskii type'; cf.\ \cite[Thm. 5.2]{DZ}), with local rate function
\[I_{x,{\bs g}}(y):=\sup_\vt\left(\vt y - \Lambda_{x,{\bs g}}(\vt)\right).\]
This concretely means that, with
$Y^\circ(t):= N^{-1}Y(t)$ and $t\in[0,T]$, and under mild regularity conditions on the set $A$,
\[\lim_{N\to\infty}\frac{1}{N} \log {\mathbb P}( Y^\circ(\cdot) \in A) =-\inf_{f\in A} {\mathbb I}_T(f),\]
with
\[{\mathbb I}_T(f):=\inf_{{\bs g}(\cdot)\in U_T}\left(\int_0^T I_{f(s),{\bs g}(s)}(f'(s)){\rm d}s+{\mathbb J}_T({\bs g})\right).\]
A formal derivation of this {\sc ldp} is beyond the scope of this paper.

\section{Erd\H{o}s-R\'enyi graphs with resampling} \label{ERR}
An alternative dynamic Erd\H{o}s-R\'enyi model (in discrete time) can be defined as follows; we refer to it as a Erd\H{o}s-R\'enyi graph with resampling.  Let the $N$ edges  alternate between two states: the edge has the value $0$ when the corresponding edge is absent and $1$ when it exists. In slot $m$, let the transition matrix of the presence of any of the $N$ edges be given by\[\left(\begin{array}{cc}P_m&1-P_m\\1-R_m&R_m\end{array}\right),\]
where the sequence $(P_m,R_m)_{m\in{\mathbb N}}$ consists of i.i.d.\ vectors in $(0,1)^2$; we note that $P_m$ and $R_m$ (for a given time $m$, that is) are {\it not} necessarily assumed independent.
It is stressed that the samples in slot $m$, i.e., $P_m$ and $R_m$, hold for {\it any} of the edges --- as a consequence, the individual edges (each of them alternating between absent and present) evolve dependently, as intended. 

In this section we find the counterparts for the resampling model of all results that we derived for the regime switching model of Section \ref{RS}. To make notation compact, let $(P,R)$ denote a generic sample of $(P_m,R_m)$.

\subsection{Generating function}

Let us now analyze the object
$\varphi_k(z) := {\mathbb E}\left(z^{Y_{m}}\,|\, Y_{m-1}=k\right).$
Realize that $Y_{m}$ is the sum of (i)~the edges that were present at time $m-1$ and still are at $m$, and (ii)~the edges that were not there at $m-1$ but do appear at $m$. Both obey a binomial distribution (with appropriately chosen parameters). More precisely,
\[ \varphi_k(z) ={\mathbb E}\left(\sum_{\ell=0}^{N-k} \binom{N-k}{\ell} (1-P_m)^\ell P_m^{N-k-\ell} z^\ell \cdot\sum_{\ell=0}^{k} \binom{k}{\ell} R_m^\ell (1-R_m)^{k-\ell} z^\ell\right),\]
which simplifies to
\[{\mathbb E}\left(\left((1-P_m)z+P_m\right)^{N-k}\cdot\left(R_mz+1-R_m\right)^k\right).\]
Now consider the stationary random variable $Y$, through its $z$-transform $\varphi(z):= 
{\mathbb E}\,z^{Y}.$ Based on the above computation, we have found the following fixed-point equation:
\begin{equation}
\label{FPE}\varphi(z) = {\mathbb E}\left(((1-P)z +P)^N \,\varphi\left(\frac{Rz+1-R}{(1-P)z+P}\right)\right).\end{equation}

\subsection{Moments}
In this subsection, we compute the mean, variance and correlation in stationarity.

\vb\vb

\noindent {\it Mean.}
Let us first compute ${\mathbb E}\,Y$, by differentiating both sides to $z$ and plugging in $z=1$.
To this end, we define
\[\psi_1(z):= ((1-P)z +P)^N,\:\:\:\psi_2(z):= \varphi\left(\frac{Rz+1-R}{(1-P)z+P}\right).\]
We first compute a number of quantities that we need in the sequel. It takes routine calcutations to conclude that
\begin{eqnarray*}
\psi_1'(z) &=& (1-P)N ((1-P)z +P)^{N-1},\\
\psi_1''(z) &=& (1-P)^2 N(N-1) ((1-P)z +P)^{N-2},\\
\psi'_2(z) &=& \frac{P+R-1}{((1-P)z+P)^2} \varphi'\left(\frac{Rz+1-R}{(1-P)z+P}\right),\end{eqnarray*}
and
\begin{eqnarray*}\psi_2''(z)&=&-2\frac{(P+R-1)(1-P)}{((1-P)z+P)^3}\varphi'\left(\frac{Rz+1-R}{(1-P)z+P}\right)
+\\
&&\frac{(P+R-1)^2}{((1-P)z+P)^4} \varphi''\left(\frac{Rz+1-R}{(1-P)z+P}\right).\end{eqnarray*}
As a consequence, 
\[\psi'_1(1) = (1-P)N,\:\:\:\psi''_1(1)=(1-P)^2N(N-1),\:\:\:\psi'_2(1)=(P+R-1)\varphi'(1),\]
and
\[
\psi_2''(1)= -2(P+R-1)(1-P)\varphi'(1)+(P+R-1)^2\varphi''(1) .\]
Regarding the first moment of $Y$, we obtain the equation 
$\alpha:=\varphi'(1) = {\mathbb E} \,\psi_1'(1) +  {\mathbb E} \,\psi_2'(1),$ or equivalently
$\alpha = N(1-{\mathbb E}\,P)+\alpha({\mathbb E}\,P+{\mathbb E}\,R-1),$
and hence
\begin{equation}
\label{MEAN}
\alpha =N \,\frac{1-{\mathbb E}\,P}{2-{\mathbb E}\,P-{\mathbb E}\,R}.\end{equation}

\vb

\noindent {\it Variance.}
We now evaluate the quantity \[\beta:={\mathbb E}\,Y(Y-1) =\varphi''(1)= {\mathbb E} \,\psi_1''(1) + 2\,
{\mathbb E} \,\psi_1'(1)\psi_2'(1)+
 {\mathbb E} \,\psi_2''(1).\] We thus obtain that $\beta$ equals
 \[N(N-1)\,{\mathbb E}\left((1-P)^2\right) + 2( N-1) \,\alpha\, {\mathbb E}\left((P+R-1)(1-P)\right)+\beta \,{\mathbb E}\left((P+R-1)^2\right),\]
 and therefore
 \[\beta = \frac{N(N-1)\,{\mathbb E}\left((1-P)^2\right) + 2( N-1) \,\alpha\, {\mathbb E}\left((P+R-1)(1-P)\right)}{1-
{\mathbb E}\left((P+R-1)^2\right) }.\]
As a consequence, ${\mathbb V}{\rm ar}\,Y$ equals
\[  \alpha-\alpha^2+\frac{N(N-1)\,{\mathbb E}\left((1-P)^2\right) + 2( N-1) \,\alpha\, {\mathbb E}\left((P+R-1)(1-P)\right)}{1-
{\mathbb E}\left((P+R-1)^2\right) }.\]
It takes an elementary but tedious computation to very that if $P$ and $R$ equal (deterministically) $p$ and $r$, respectively, then this variance reduces to $N\pi_0\pi_1$, as desired.

We also conclude that ${\mathbb V}{\rm ar}\,Y$ grows essentially quadratic in $N$. Indeed, it follows by standard computations that, with $\bar P:=1-P$ and $\bar R:=1-R$, 
\begin{equation}
\label{VAR}{\mathbb V}{\rm ar}\,Y = \gamma_1 N^2 +\gamma_2 N,\end{equation}
where
\[\gamma_1 = \frac{{\mathbb E} (\bar R^2) ({\mathbb E}\,\bar P)^2 -
2 \,{\mathbb E} (\bar P\bar R) {\mathbb E}\,\bar P\,{\mathbb E}\,\bar R
+  {\mathbb E} (\bar P^2) ({\mathbb E}\,\bar R)^2 }{\left(1-
{\mathbb E}\left((\bar P+\bar R-1)^2\right)\right) \left({\mathbb E}\,\bar P+{\mathbb E}\,\bar R\right)^2},\]
and
\[\gamma_2 = \frac{-\,{\mathbb E} (\bar R^2) \,{\mathbb E}\,\bar P
+ 2 \,{\mathbb E}\,\bar P\,{\mathbb E}\,\bar R
-{\mathbb E} (\bar P^2) \,{\mathbb E}\,\bar R
}{\left(1-
{\mathbb E}\left((\bar P+\bar R-1)^2\right)\right) \left({\mathbb E}\,\bar P+{\mathbb E}\,\bar R\right)}.\]
Notice that $\gamma_1$ and $\gamma_2$ are symmetric in $P$ and $R$, as desired, and observe that $\gamma_1\ge 0$ (with equality only if $P$ and $R$ are deterministic).
We conclude that no standard CLT applies (which would require that ${\mathbb V}{\rm ar}\,Y$ grows linearly in $N$) unless $P$ and $R$ are deterministic.

\vb

\noindent {\it Correlation.}
We now focus on computing the limit of covariance ${\mathbb C}{\rm ov} (Y_m,Y_{m+1})$ as $m\to\infty.$
\iffalse
A way to go about this is by analyzing the object
\[\check\varphi_k(w,z) := {\mathbb E}\left(w^{Y_{m-1}} \,z^{Y_{m}}\,|\, Y_{m-2}=k\right).\]
With 
\[\xi_{k,\ell}(w,r):= \binom{k}{\ell}  r^\ell (1-r)^{k-\ell}\,w^\ell,\]
it takes an elementary conditioning argument to conclude   that
\begin{eqnarray*}\lefteqn{\hspace{-1mm}\xi_{k,\ell}(w,r) = {\mathbb E}\left[  \left(\sum_{\ell=0}^k \xi_k(w,R_{m-1})\left(\sum_{\ell'=0}^\ell \xi_{\ell,\ell'}(z,R_m)\right)\left(\sum_{\ell'=0}^{k-\ell}\xi_{k-\ell,\ell'}(z,1-P_m) \right)\right)\right.}\\
&&\hspace{1.3cm}\left.  \left(\sum_{\ell=0}^{N-k} \xi_{N-k}(w,1-P_{m-1})\left(\sum_{\ell'=0}^\ell \xi_{\ell,\ell'}(z,R_m)\right)\left(\sum_{\ell'=0}^{N-k-\ell}\xi_{N-k-\ell,\ell'}(z,1-P_m) \right)\right)\right].
\end{eqnarray*}
This simplifies to 
\[{\mathbb E}\left[\left(
{\mathscr A}(w,z)
\right)^k \left({\mathscr B}(w,z)\right)^{N-k}\right],\]with 
\begin{eqnarray*}
{\mathscr A}(w,z)&:=& R_{m-1}w(R_mz+1-R_m)+(1-R_{m-1})((1-P_m)z+P_m),\\
{\mathscr B}(w,z)&:=& (1-P_{m-1})w(R_mz+1-R_m)+P_{m-1}((1-P_m)z+P_m).
\end{eqnarray*}
Conclude that, in stationarity, with the $\varphi(\cdot)$ we characterized above,
\[\check\varphi(w,z) := \lim_{m\to\infty}
{\mathbb E}\left(w^{Y_{m-1}} \,z^{Y_{m}} \right) = {\mathbb E}\left(
\left({\mathscr B}(w,z)\right)^N \varphi\left(\frac{{\mathscr A}(w,z)}{{\mathscr B}(w,z)}\right)\right).
\]
This provides us with the full  distribution corresponding to two subsequent time epochs in stationarity. If we are just interested in the covariance, a simpler approach is possible.\fi
Observe that
\[\lim_{m\to\infty} {\mathbb C}{\rm ov} (Y_m,Y_{m+1}) = 
\lim_{m\to\infty} \sum_{k=0}^N k \,{\mathbb E}(Y_{m+1}\,|\, Y_m = k) \,{\mathbb P}(Y_m = k) - ({\mathbb E}\,Y)^2,\]
which, in self-evident notation, reads
\[\sum_{k=0}^N k \,{\mathbb E}({\mathbb B}{\rm in} (k,R))\,{\mathbb P}(Y= k) 
+\sum_{k=0}^N k\,{\mathbb E}({\mathbb B}{\rm in} (N-k,1-P) )  \,{\mathbb P}(Y= k) - ({\mathbb E}\,Y)^2.\]
This reduces to
\[\,{\mathbb E}R\,\sum_{k=0}^N k^2 \,{\mathbb P}(Y= k) 
+(1-\,{\mathbb E}P)\sum_{k=0}^N k(N-k)  \,{\mathbb P}(Y= k) - ({\mathbb E}\,Y)^2,\]
so that we obtain
\[\lim_{m\to\infty} {\mathbb C}{\rm ov} (Y_m,Y_{m+1}) =({\mathbb E}P+{\mathbb E}R-1) {\mathbb E}(Y^2) + (1-\,{\mathbb E}P)N \,{\mathbb E}\,Y - ({\mathbb E}\,Y)^2,\]
which we can evaluate from the expressions for ${\mathbb E}\,Y$ and ${\mathbb V}{\rm ar}\,Y.$ 

\subsection{Diffusion results under scaling}
We now consider the following scaling: for some $\delta>0$ we put
\begin{equation}\label{PR}
P= 1-\eta/N^\delta,\:\:\:\:R= 1-\zeta/N^\delta,\end{equation}  where $\eta$ and $\zeta$ are non-negative random variables. 
The resulting model has some built-in `inertia': for $N$ large, the process has the inclination to stay in the same configuration. 
The mean number of vertices is
$N\,\bar\varrho,$ with\[\bar\varrho:=\frac{{\mathbb E}\,\eta}{{\mathbb E}\,\eta+{\mathbb E}\,\zeta},\]
irrespective of the value of $\delta$.
When analyzing the variance, however, the revealing issue is that the value of $\delta$ has crucial impact. More specifically, a minor computation tells us that ${\mathbb V}{\rm ar}\,Y$ essentially  reads
\[N\,\bar\varrho\,(1-\,\bar\varrho\,)+N^{2-\delta} \frac{{\mathbb E} (\zeta^2) ({\mathbb E}\,\eta)^2 -
2 \,{\mathbb E} (\eta\zeta) {\mathbb E}\,\eta\,{\mathbb E}\,\zeta
+  {\mathbb E} (\eta^2) ({\mathbb E}\,\zeta)^2 }{2({\mathbb E}\,\eta+{\mathbb E}\,\zeta)^3}  .\]
Note that, due to the inertia that we incorporated,  the variance is smaller than in the unscaled model, where the variance was effectively proportional to $N^2$. 
Observe from the above expression that there is  a dichotomy that resembles the one we came across in Section \ref{RS}, with  some sort of transition at $\delta = 1.$ For $\delta>1$ 
the standard deviation scales as $\sqrt{N}$, whereas for  $\delta<1$ it scales as $N^{1-\delta/2}$. An intuitive explanation is that in the regime of relatively few transitions (i.e., $\delta>1$) the system's inertia is so strong that its steady-state essentially behaves as an Erd\H{o}s-R\'enyi graph with the probability that an edge exists being given by
$\bar\varrho$.
In the regime with relatively many transitions (i.e., $\delta<1$), on the contrary, the (co-)variances play a role, in the sense that the increased variability caused by the resampling has impact; the limiting object is not of Erd\H{o}s-R\'enyi-type.

Along the same lines, an elementary computation yields that the covariance between the numbers of edges at two subsequent epochs (in stationarity) behaves as 
\[{\mathbb V}{\rm ar}\,Y\left(1-\frac{{\mathbb E}\,\eta+{\mathbb E}\,\zeta}{N^\delta}\right);\]this correlation coefficient essentially reads 
$1-({{\mathbb E}\,\eta+{\mathbb E}\,\zeta}){N^{-\delta}}$ (for $N$ large).

\vb

\noindent
{\it A related continuous-time model.}
In the remainder of this subsection we consider a specific explicit continuous-time model in which we can embed the discrete-time model discussed above, and in particular the scaling (\ref{PR}). To this end, we first describe the model without scaling, and then include the scaling.
 
Let, at time $s$, $M(s)\geqslant 0$ be the hazard rate of an existing vertex becoming inactive; likewise, $\Lambda(s)\geqslant 0$ is the hazard rate corresponding with a non-existing vertex becoming active. Here $M(s)$ and $\Lambda(s)$ are piecewise constant stochastic processes: for some $\Delta>0$,
\[\Lambda(s) = \Lambda_i \,1_{\{(i-1)\Delta \leqslant s  < i \Delta\}},\:\:\:
M(s) =M_i \,1_{\{(i-1)\Delta \leqslant s  < i \Delta\}},\]
where $(M_i,\Lambda_i)_{i\in{\mathbb N}}$ is a sequence of i.i.d.\ bivariate random vectors such that both $ {\mathbb V}{\rm ar}\, \Lambda$ and $ {\mathbb V}{\rm ar}\, M$ are finite. Let $Y(t)$ be the number of vertices at time $t$, and $Y$ its stationary counterpart. As it turns out, we can reuse quite a few results from the previous subsections, using the identification $Y(m\Delta) = Y_m.$
In particular, it is seen that $\varphi(z) :={\mathbb E}\, z^Y$ satisfies (\ref{FPE}), with 
\[P := \frac{M}{\Lambda+M}+\frac{\Lambda}{\Lambda+M}e^{-(\Lambda+M)\Delta} ,\:\:R:=\frac{\Lambda}{\Lambda+M}+\frac{M}{\Lambda+M}e^{-(\Lambda+M)\Delta}.\]
We thus obtain from (\ref{MEAN})
\[{\mathbb E}\,Y = N\, {\mathbb E}\left(\frac{\Lambda}{\Lambda+M}\left(1-e^{-(\Lambda+M)\Delta}\right)\right)\left/{\mathbb E}\left(1-e^{-(\Lambda+M)\Delta}\right)\right.\]
%The limiting values for $\Delta\downarrow 0$ and $\Delta\to 0$ have an appealing intuitive explanation:
%\[\lim_{\Delta\downarrow 0} {\mathbb E}\,Y = \frac{{\mathbb E}\,\Lambda}{{\mathbb E}\,\Lambda+{\mathbb E}\,M},\:\:\:\:\:\:\lim_{\Delta\to\infty}{\mathbb E}\,Y ={\mathbb E}\left(\frac{\Lambda}{\Lambda+M}\right).\]
Similarly, we can compute the variance by (\ref{VAR}). 

Now we describe how to scale  this model. The idea  is to scale  $\Delta\mapsto 1/N^{\delta}$, and to consider the regime in which we let $N$ grow large, i.e., the transition rates are frequently resampled (and simultaneously the number of potential edges $N$ grows). It is immediate that $P$ and $R$ fulfill (\ref{PR}) with $\eta=\Lambda$ and $\zeta=M.$ We obtain that 
${\mathbb E}\,Y$ tends to $\bar\varrho :={\mathbb E}\,\Lambda/{\mathbb E}\,\Gamma$, where  $\Gamma:=\Lambda+M$.
In addition, $ {\mathbb V}{\rm ar}\, Y$ satisfies the expansion $N\,\bar\varrho\,(1-\,\bar\varrho\,)+N^{2-\delta} v+o(N^{\max\{1,2-\delta\}})$, where  
\begin{eqnarray*}v&:=&\frac{1}{2\,{\mathbb E}\,\Gamma} \left(\,\bar\varrho\,^2\, {\mathbb V}{\rm ar}\, M - 2\,\bar\varrho\,(1-\,\bar\varrho\,)\, {\mathbb C}{\rm ov}\,(\Lambda,M) +
(1-\,\bar\varrho\,)^2\, {\mathbb V}{\rm ar}\, \Lambda\right)\\&=&\frac{1}{2\,{\mathbb E}\,\Gamma} {\mathbb V}{\rm ar}\left(\Lambda-\,\bar\varrho\,\Gamma\right).\end{eqnarray*}
The proof of a functional central limit theorem is very similar to the one for the regime switching model in Section \ref{RS}; we therefore restrict ourselves to the key steps. With $P_1(\cdot)$ and $P_2(\cdot)$ as before,
\[Y(t) = P_1\left(\int_0^t \Lambda(s) (N-Y(s)){\rm d}s\right) - P_2\left(\int_0^t M(s) Y(s) {\rm d}s\right),\]
so that, for some martingale $K(t)$,
\[{\rm d}Y(t) =\Lambda(t) (N-Y(t)){\rm d}t-M(t) Y(t){\rm d}t +{\rm d} K(t).\]
Then $\bar Y(t)$ is defined as in (\ref{barY}), with  $\varrho(t):=\bar\varrho\cdot (1-\exp(-t\,{\mathbb E}\,\Gamma)).$ We define, with $\Gamma(s) =\Lambda(s)+M(s)$,
\[W(t):=e^{\Gamma_+(t)}\bar Y(t),\:\:\:\:\mbox{with}\:\:\:
\Gamma_+(t):=\int_0^t \Gamma(s){\rm d}s .\]
After a few steps, this leads to the stochastic differential equation, 
\[{\rm d}W(t) = 
e^{\Gamma_+(t)} \left(\sqrt{N}\left((\Lambda(t)-{\mathbb E}\,\Lambda)-
\varrho(t)(\Gamma(t)-{\mathbb E}\,\Gamma)\right){\rm d}t+\frac{{\rm d}K(t)}{\sqrt{N}}\right).\]
Consider the two terms in the previous display. For the first term, as $N\to\infty$,
\[\int_0^\cdot \sqrt{N} e^{ \Gamma_+(s)}
\big((\Lambda(s)-{\mathbb E}\,\Lambda)-
\varrho(s)(\Gamma(s)-{\mathbb E}\,\Gamma)\big)
 {\rm d}s\to \int_0^\cdot e^{s\,{\mathbb E}\,\Gamma }{\rm d} G(s),\]
where $G(\cdot)$ satisfies
\begin{equation}
\label{defg}\langle G\rangle_t = g(t):=\int_0^t {\mathbb V}{\rm ar}\, (\Lambda -\varrho(s)\Gamma) {\rm d}s;
\end{equation}
to see this note that, almost surely, uniformly on compacts, as $N\to\infty$,
\[ e^{ \Gamma_+(s)}=  \exp\left(\frac{1}{N}\sum_{i=1}^{sN} (\Lambda_i+M_i)\right) \to
\exp\left(s\,{\mathbb E}\,\Gamma\right),\]
and use this in combination with the (classical) functional central limit theorem for the random walk with i.i.d.\ increments \cite[Thm. 4.3.5]{WHITT}.
For the second term, as $N\to\infty$, due to the definition of the martingale $K(\cdot)$,
\[\int_0^\cdot \frac{1}{\sqrt{N}}e^{ \Gamma_+(s)} {\rm d}K(s)  \to 
 \int_0^\cdot e^{\gamma\s s}{\rm d} H(s),\]
where $H(\cdot)$ is such that
\begin{equation}
\label{defh}\langle H\rangle_t = h(t):={\mathbb E}\,\Lambda \int_0^t (1-\varrho(s)){\rm d}s +{\mathbb E}\,M \int_0^t \varrho(s){\rm d}s.\end{equation}
Combining the two terms studied above, it thus follows that, as $N\to\infty$, $W(\cdot)$ weakly converges to $W_\infty(\cdot)$, which is the solution to the stochastic differential equation (\ref{defW}), but now with the
$g(\cdot)$ and $h(\cdot)$ given by (\ref{defg}) and (\ref{defh}), respectively.  
We obtain the following result.

\begin{theorem}
$\bar Y(\cdot)$  converges weakly to $\bar Y_\infty(\cdot)$, which is the solution to the stochastic differential equation $(\ref{limde})$,
with 
$g(\cdot)$ and $h(\cdot)$ given by $(\ref{defg})$ and $(\ref{defh})$, respectively.
\end{theorem}

\begin{remark}
For large  $t$ (`in stationarity'), this stochastic differential equation essentially behaves as
\[ {\rm d}\bar Y_\infty(t)= - {\mathbb E}\,\Gamma\cdot\bar Y_\infty(t)\,{\rm d}t + \sqrt{2\,{\mathbb E}\,\Gamma\cdot\bar\varrho(1-\bar\varrho)+2\,{\mathbb E}\,\Gamma\cdot v}\,{\rm d}B(t),\]
corresponding with an {\sc ou} process with mean $0$ and variance $\bar\varrho\,(1-\bar\varrho) + v$. Note that this is in line  with what we found, plugging in $\delta=1$, in the expansion $N\,\bar\varrho\,(1-\,\bar\varrho\,)+N^{2-\delta} v+o(N^{\max\{1,2-\delta\}})$. Regarding the cases $\delta<1$ and $\delta>1$ a reasoning similar to that in Remark \ref{rem2} applies. 
\end{remark}

\subsection{Large deviations results under scaling}

\iffalse Above we observed that we have the following recursion:
\[Y_{m} = Y_{m-1} - B_1+B_2,\:\:\:\mbox{with}\:\:\:
B_1:={\mathbb B}{\rm in} (Y_{m-1},1-R_m),\:\:B_2:= {\mathbb B}{\rm in} (N-Y_{m-1},1-P_m).\]
From now on we consider the case $\delta = 1$ in more detail.
Denoting $Y_{m-1}=Nx$, we obtain the approximation
\[ Y_{m} - Y_{m-1}= \mu + U\sigma,\]
with $U$ being standard Normal, and (by conditioning on the values of $P_m$ and $R_m$)
\[\mu:= (1-x){\mathbb E}\,\eta -x\,{\mathbb E}\,\zeta,\]
\[\sigma^2 := {(1-x)^2{\mathbb V}{\rm ar}\, \eta + 2x(1-x) {\mathbb C}{\rm ov}( \eta,\zeta) +x^2{\mathbb V}{\rm ar}\, \zeta  +(1-x){\mathbb E}\,\eta +x\,{\mathbb E}\,\zeta }.\]
Challenge ii:} Derive an fclt} for this Markov chain, or just a clt} for its stationary distribution. My gut feeling is that an fclt} with an ou} limit could come out. \fi

\vb

The above computations focused on the mean, variance, and correlation under the scaling proposed. 
We now consider rare events.
Another straightforward calculation yields for the cumulant function, assuming $Nx$ to be integer,
\begin{eqnarray*}\lefteqn{\log {\mathbb E}\,\exp\left(\vt(Y_{m} - Y_{m-1})\,|\,Y_{m-1} = Nx\right) }\\&=&\log {\mathbb E}\left(\left(e^{-\vt}(1-R_m)+R_m)\right)^{Nx}
\left(e^{\vt}(1-P_m)+P_m)\right)^{N(1-x)}\right),\end{eqnarray*}
which, for $\delta = 1$, converges to 
\begin{eqnarray*}\Lambda_x(\vt)&:=&\log {\mathbb E}\exp\left({x\zeta (e^{-\vt} -1) +(1-x)\eta (e^\vt-1)} \right) \\&=&\log M\left( x (e^{-\vt} -1),(1-x)(e^\vt-1) \right)
%=\log {\mathbb E}\left(e^{x\zeta e^{-\vt}+(1-x)\eta e^\vt} \right)-1
,\end{eqnarray*}
where $M(\cdot,\cdot)$ is the joint moment generating function of the random variables $\zeta$ and $\eta$ (assuming that it exists).
One thus finds a sample-path {\sc ldp} where the local rate function is given by
\[I_x(y):=\sup_\vt\left (\vt y -\Lambda_x(\vt) \right).\]
More precisely, with
$Y^\circ(t):= N^{-1}Y_{\lfloor Nt\rfloor}$ and $t\in[0,T]$, and under mild regularity conditions on the set $A$,
\[\lim_{N\to\infty}\frac{1}{N} \log {\mathbb P}( Y^\circ(\cdot) \in A) =-\inf_{f\in A} {\mathbb I}_T(f),\:\:\:\mbox{with}\:\:\:{\mathbb I}_T(f):=\int_0^T I_{f(s)}(f'(s)){\rm d}s.\]

\vspace{-0.6cm}
\begin{figure}
\includegraphics[width=7.8cm]{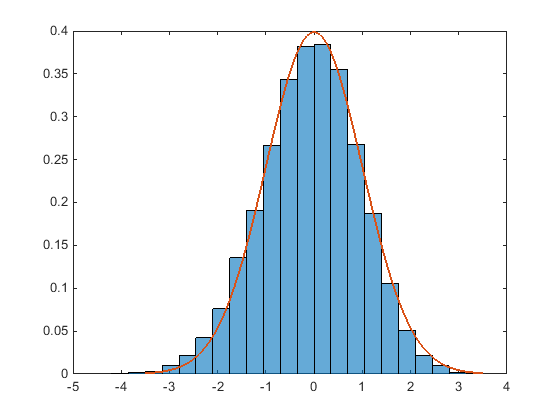}
\includegraphics[width=7.8cm]{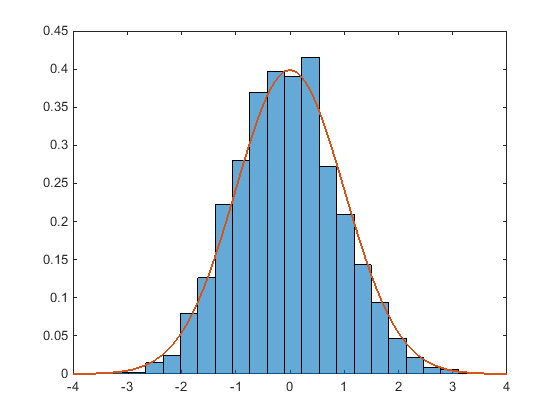}
\caption{\label{Fig1}Left panel: histogram of $\bar Y$ for situation (A). Right panel:  histogram of $\bar Y$ for situation (B). In both cases we took $N=45$.}
\end{figure}

\vspace{-0.6cm}

\section{Numerical illustration}\label{NUM}
In this section we include a number of illustrative examples that assess the applicability of the diffusion limits.
We consider two situations; in both cases we take $\delta=1$. (A)~In the first situation we consider the regime switching model of Section \ref{RS}. The background process has two states, with $q_{12}=2$ and $q_{21}=3$; in addition $\lambda_1=0.3,$ $\lambda_2=0.5$, $\mu_1=1$, and $\mu_2=0.1$. Using the formulae we derived in Section \ref{RS}, we find ${\mathbb E}\,Y=0.762\,N$ and ${\mathbb V}{\rm ar}\,Y= 0.182\,N$.
(B) The second situation corresponds to the resampling model of Section \ref{ERR}. More, specifically, $M$ has a uniform distribution on $[0,3]$ and $\Lambda$ a uniform distribution on $[0,5]$. It is readily checked that ${\mathbb E}\,Y= 0.625\,N$ and ${\mathbb V}{\rm ar}\, Y= 0.308\,N.$

In Fig.\ 1 histograms are presented for the random variable
\[\bar Y:= \frac{Y- {\mathbb E}\,Y}{\sqrt{{\mathbb V}{\rm ar}\,Y}}.\]
The number of experiments the estimates are based upon equals the number of this {\sc lncs} volume. Each simulation experiment starts with an empty system, and is then run for a sufficiently long time such that the process has reached equilibrium. The red curves in Fig.\ 1 correspond to the density of the standard Normal distribution. The figures confirm the convergence to the Normal distribution. 

In Fig.\ 2 typical sample paths are depicted, illustrating the {\sc ou}-like mean-reverting behavior. The red curves correspond to the mean of $Y(t)$. 

\begin{figure}
\includegraphics[width=7.8cm]{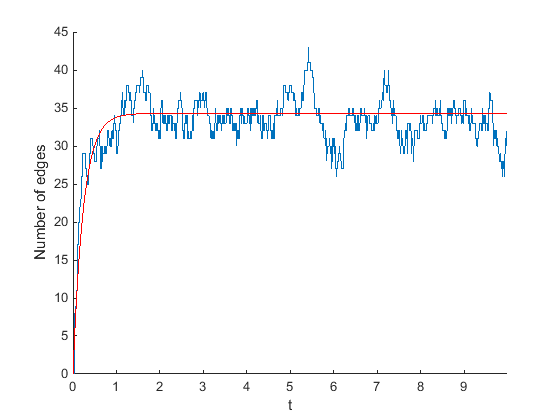}
\includegraphics[width=7.8cm]{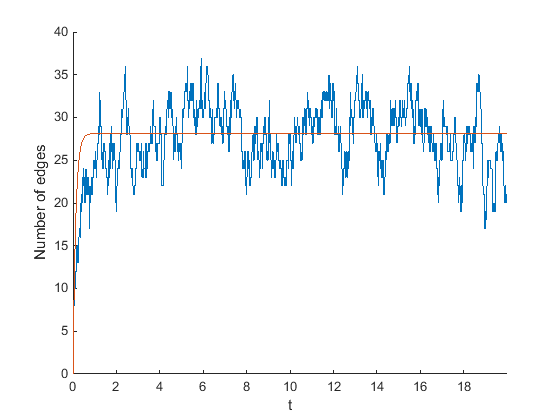}
\caption{\label{Fig1}Left panel: sample path of $Y(\cdot)$ for situation (A). Right panel:  sample path of $Y(\cdot)$ for situation (B). In both cases we took $N=45$.}
\end{figure}

\vspace{-0.2cm}

\section{Discussion and concluding remarks}
In this paper we have discussed distributional properties of the number of edges in a dynamic Erd\H{o}s-R\'enyi graph. 
We have considered two variants: one with the underlying mechanism being based on regime switching, and the other in which the transition probabilities are resampled at equidistant points in time.  For both models we have succeeded in obtaining fairly explicit results for various transient and stationary quantities. Under a specific scaling a functional central limit theorem was established. 

There is an interesting relation between the models considered in this paper and two-node closed queueing networks. In such closed networks a fixed number of jobs, say $N$, move between an active state (`in service') and an inactive state (`waiting'). Such models (but without regime switching or resampling) have been intensively studied in the literature in the context of so-called Engset models \cite{ENG}; see e.g.\ \cite{CH} and references therein. 

Topics for future research may relate to other graph metrics than the total number of edges. In the introduction, we mentioned that \cite{TYA} considers the behavior of the Betti number, but one could also think of e.g.\ the evolution of the number of wedges or triangles in the random graph. In addition, one may wonder  under what conditions the dynamic random graph in which the edges (independently) alternate between present and absent is almost surely connected; one would expect that if this alternating process is `sufficiently fast' and the stationary up-probability is larger than $\log n/n$, this should be the case.

\vb

\noindent {\small {\bf Acknowledgment} -- The authors thank Frank den Hollander (Leiden) for useful discussions.}

%{\small \section*{Acknowledgment}

%
% ---- Bibliography ----
%

\section*{Appendix}
We now prove Lemma \ref{lalemma}. We do so by establishing the claim for $k=1$; plugging in $k\,\Gamma$ for $\Gamma$ yields the stated.  Write $F_\infty:=(\gamma\s)^{-1}{\bs 1}{\bs \pi}^{\rm T}$ and abbreviate $F_N:=F_{N,1}.$

As $Q$ has a kernel of dimension 1, we can factorize $Q$ as $Q=AB$, where $A\in\rr^{d\times(d-1)}$ is of full column rank and $B\in\rr^{(d-1)\times d}$ is of full row rank. It is not hard  to show that 
%in the present situation 
$BA$ is an invertible matrix. Moreover, every element in the right kernel of $Q$ is a multiple of ${\bs 1}$ and, likewise, every element in the left kernel of $Q$ is a multiple of ${\bs \pi}^{\rm T}$. 

Applying the Sherman-Morrison formula to $F_N=(\Gamma-NAB)^{-1}$, we find
\begin{equation}\label{eq:SM}
F_N=\Gamma^{-1}+\Gamma^{-1}A\left(\frac{I_{d-1}}{N}-B\Gamma^{-1}A\right)^{-1}B\Gamma^{-1}.
\end{equation}
Taking the limit for $N\to\infty$, we arrive at 
\begin{equation}\label{eq:finfty}
F_\infty=\Gamma^{-1}-\Gamma^{-1}A\,(B\Gamma^{-1}A)^{-1}B\Gamma^{-1}, 
\end{equation}
where the invertibility of $B\Gamma^{-1}A$ is due to ${\bs \gamma}\s>0$.  One sees that  $F_\infty A=0$ and $BF_\infty=0$. Hence $F_\infty$ belongs to the left kernel of $A$ and to the right kernel of $B$, so $F_\infty=c\,{\bs 1}{\bs \pi}^{\rm T}$ for some $c\in\rr$. One also has $F_\infty\Gamma{\bs 1}=F_\infty{\bs \gamma}={\bs 1}$, and hence $c=({\gamma\s})^{-1}$, which gives the desired result for $\lim_{N\to\infty}F_N$.

We proceed by proving the expansion. Inserting 
\begin{align*}
\left(\frac{I}{N}-B\Gamma^{-1}A\right)^{-1}
& = -(B\Gamma^{-1}A)^{-1}-\frac{1}{N}(B\Gamma^{-1}A)^{-2}+O(N^{-2})
\end{align*}
into~\eqref{eq:SM}, one obtains
\[
F_N=F_\infty-\frac{1}{N}\Gamma^{-1}A(B\Gamma^{-1}A)^{-2}B\Gamma^{-1}+O(N^{-2}).
\]
Let $A^+$ ($B^+$, resp.) denote any left (right, resp.) inverse of $A$ ($B$, resp.), so that $A^+A=BB^+=I_{d-1}$. Then 
\[
\Gamma^{-1}A(B\Gamma^{-1}A)^{-2}B\Gamma^{-1}= \Gamma^{-1}A(B\Gamma^{-1}A)^{-1}BB^+A^+A(B\Gamma^{-1}A)^{-1}B\Gamma^{-1}.
\]
Now, it follows from \eqref{eq:finfty} that $\Gamma^{-1}A(B\Gamma^{-1}A)^{-1}B=I-F_\infty\Gamma$ and in addition
$A(B\Gamma^{-1}A)^{-1}B\Gamma^{-1}=I-\Gamma F_\infty$. Hence,
\begin{equation}\label{eq:fn}
F_N=F_\infty-\frac{1}{N}(I-F_\infty\Gamma)B^+A^+(I-\Gamma F_\infty)+O(N^{-2}).
\end{equation}
We specialize to judicious choices of $A^+$ and $B^+$, namely
\[
A^+=
\begin{pmatrix}
I_{d-1} & 0
\end{pmatrix}
A_1^{-1},\:\: B^+=
B_1^{-1}
\begin{pmatrix}
I_{d-1} \\ 0
\end{pmatrix}\hspace{-0.6mm},\:\:\hspace{-0.6mm}\mbox{where}
\:\:\hspace{-0.6mm}A_1:=\begin{pmatrix}
A & {\bs 1}
\end{pmatrix}\hspace{-0.6mm},\:\:B_1:=\begin{pmatrix}
B \\
-{\bs \pi}^{\rm T}
\end{pmatrix}\hspace{-0.6mm},\] and $0$ stands here as well as below for a zero matrix or vector of appropriate dimensions. Both $A_1$ and $B_1$ are invertible, as an immediate consequence of the relation
$
A_1B_1  =Q-{\bs 1}{\bs \pi}^{\rm T}=-(D+{\bs 1}{\bs \pi}^{\rm T})^{-1}.$ In addition,
\[B^+A^+  =B_1^{-1}\begin{pmatrix}
I_{d-1} & 0 \\ 0 & 0
\end{pmatrix}A_1^{-1},\:\:
B_1{\bs 1} = -\begin{pmatrix}
0 \\ 1
\end{pmatrix},\:\:
{\bs \pi}^{\rm T} A_1  = \begin{pmatrix}
0 & 1
\end{pmatrix}.
\]
A straightforward computation gives with the above relations
\[
B^+A^+  =B_1^{-1}A_1^{-1}-B_1^{-1}
\begin{pmatrix}
1 \\
0
\end{pmatrix}
\begin{pmatrix}
1 & 0
\end{pmatrix}
A_1^{-1} = -(D+{\bs 1}{\bs \pi}^{\rm T})+{\bs 1}{\bs \pi}^{\rm T} = -D.\]The result (for $k=1$) now follows from \eqref{eq:fn}.  
\hfill$\Box$

\end{document}